\documentclass[10pt]{amsart}

\let\hat=\widehat

\newcommand{\Fq}{{\mathbf F}_q}
\newcommand{\Fqq}{{\mathbf F}_{q^2}}
\newcommand{\BC}{{\mathbf C}}
\newcommand{\BF}{{\mathbf F}}

\newcommand{\BQ}{{\mathbf Q}}
\newcommand{\BR}{{\mathbf R}}
\newcommand{\BZ}{{\mathbf Z}}
\newcommand{\BFqb}{{\overline{\BF}_q}}

\newcommand{\CC}{{\mathcal C}}
\newcommand{\CO}{{\mathcal O}}

\newcommand{\Aid}{{\mathfrak A}}
\newcommand{\Bid}{{\mathfrak B}}
\newcommand{\Cid}{{\mathfrak C}}
\newcommand{\Did}{{\mathfrak D}}
\newcommand{\Eid}{{\mathfrak E}}
\newcommand{\fid}{{\mathfrak f}}
\newcommand{\pid}{{\mathfrak p}}

\newcommand{\eps}{\varepsilon}
\newcommand{\Dp}{D^{+}}
\newcommand{\Gp}{G^{+}}
\newcommand{\Kp}{K^{+}}

\newcommand{\Rp}{R^{+}}
\newcommand{\Up}{U^{+}}
\newcommand{\fidp}{\fid^{+}}
\newcommand{\pibar}{\overline{\pi}}
\newcommand{\Aidbar}{\overline{\Aid}}
\newcommand{\Bidbar}{\overline{\Bid}}

\DeclareMathOperator{\Aut}{Aut}
\DeclareMathOperator{\End}{End}
\DeclareMathOperator{\Gal}{Gal}
\DeclareMathOperator{\Cl}{Cl}
\DeclareMathOperator{\Clp}{Cl^+}
\DeclareMathOperator{\Tor}{Tor}

\newcommand\ra{\rightarrow}
\newcommand\lra{\longrightarrow}
\newcommand\da{\downarrow}
\newcommand\bda{\big\da}

\newtheorem{theorem}{Theorem}
\newtheorem{lemma}[theorem]{Lemma}

\newtheorem{proposition}[theorem]{Proposition}

\theoremstyle{remark}
\newtheorem*{rem}{Remark}
\newtheorem*{ack}{Acknowledgment}

\begin{document}

\title[The nonexistence of certain curves]
      {On the nonexistence of\\ certain curves of genus two}
\author{Everett W.~Howe}
\address{Center for Communications Research, 
         4320 Westerra Court, 
         San Diego, CA 92121-1967, USA.}
\email{however@alumni.caltech.edu}
\urladdr{http://alumni.caltech.edu/\~{}however/}
\date{30 January 2002}


\renewcommand{\subjclassname}{\textup{2000} Mathematics Subject Classification}


\begin{abstract}
We prove that if $q$ is a power of an odd prime then
there is no genus-$2$ curve over $\Fq$ whose Jacobian has characteristic
polynomial of Frobenius equal to $x^4 + (2-2q)x^2 + q^2$.
Our proof uses the Brauer relations in a biquadratic extension of $\BQ$
to show that every principally polarized abelian surface over $\Fq$
with the given characteristic polynomial splits over $\Fqq$ as a product
of polarized elliptic curves.
\end{abstract}

\maketitle

\section{Introduction}
\label{S-intro}
Recently, Maisner and Nart~\cite{maisner-nart} compiled data on every curve
of genus $2$ over every finite field of cardinality at most~$41$.
One of the facts they gleaned from analyzing the results of their computation
is that for every odd prime power $q \le 41$,
there is no genus-$2$ curve over $\Fq$ that has exactly $q + 1$ points over $\Fq$
and exactly $q^2 - 4q + 5$ points over $\BF_{q^2}$.
With further computation they verified that the same statement holds
for all odd prime powers $q\le 61$, but they were unable
to use any of the usual arguments (mentioned below) to explain this 
observed fact.
The purpose of this paper is to provide a new argument that 
explains Maisner and Nart's observation and that shows that it
holds for all odd prime powers~$q$.

Instead of speaking about numbers of rational points, 
Maisner and Nart expressed their results in terms of the
characteristic polynomials of the curves they enumerated. 
(By the {\em characteristic polynomial\/} of a curve we mean the characteristic
polynomial of the Frobenius endomorphism of the Jacobian of
the curve.)  If for every prime power $q$ we let 
$f_q$ denote the polynomial $x^4 + (2-2q)x^2 + q^2$, 
then Maisner and Nart's observation is equivalent to the statement that
for every odd prime power $q\le 61$,
no genus-$2$ curve over $\Fq$ has characteristic polynomial~$f_q$.
Thus, our main theorem is phrased as follows:

\begin{theorem}
\label{T-main}
Let $q$ be a power of an odd prime.
There is no curve of genus $2$ over $\Fq$ whose characteristic
polynomial is equal to~$f_q$.
\end{theorem}

We should mention that when $q$ is a power of $2$ it is still true
that there is no genus-$2$ curve over $\Fq$ with characteristic polynomial~$f_q$,
but for a trivial reason: Honda-Tate theory~\cite{tate}
shows that if $q>2$ is a power of $2$ there is not even an abelian surface
with characteristic polynomial~$f_q$, let alone a Jacobian.
The remaining case, $q=2$, can easily be eliminated by
noting that no curve can have exactly $3$ points over $\BF_2$ and
exactly $1$ point over~$\BF_4$.

Up until now, there seem to have been four methods of showing that 
the characteristic polynomial of an abelian surface does 
not occur as the characteristic polynomial of a genus-$2$ curve.
Two of the methods deal directly with properties of curves:
\begin{enumerate}
\item[(1)]
For certain very small $q$, some characteristic polynomials can be excluded
because they would require the associated curve to have a negative number of points,
or fewer points over some extension field than over the base field.
(See for example~\cite[p.~353]{ruck} and the argument we gave above for the
polynomial~$f_2$.)
\item[(2)]
The polynomial $x^4 + (1-2q)x^2 + q^2$ can be excluded because 
a curve with that characteristic polynomial would have to have an
automorphism whose existence is incompatible with the number of rational
points on the curve over~$\Fqq$.
(See~\cite{howe:em}.)
\end{enumerate}
The other two methods make use of the close relationship between Jacobians
of genus-$2$ curves and principally polarized abelian surfaces.  The Jacobian 
of a genus-$2$ curve over $\BF_q$ is a principally polarized abelian
surface over $\BF_q$, and the converse is almost true: an abelian surface
over $\BF_q$ with a principal polarization is either the Jacobian
of a genus-$2$ curve over $\BF_q$ or else is isomorphic (over
the algebraic closure $\BFqb$ of $\BF_q$) to a product of polarized elliptic curves
over $\BFqb$.  So to show that there is no Jacobian corresponding to a 
given characteristic polynomial, one can show either that no principally polarized
abelian surface has that characteristic polynomial, or that 
every principally polarized abelian surface with that characteristic polynomial
is geometrically split.  This leads us to the other two methods:
\begin{enumerate}
\item[(3)]
Certain polynomials of the form $x^4 + ax^3 + (a^2 - q)x^2 + aqx + q^2$
can be excluded because one can show that the associated abelian surfaces
do not have principal polarizations.
(See~\cite[\S 13]{howe:tams}.)
\item[(4)]
Certain polynomials of the form $x^4 + (2c+1)x^3 + (c^2+c+2q)x^2 + (2c+1)qx + q^2$
that are associated with principally polarized abelian
surfaces can be excluded because one can show that the endomorphism rings
of the surfaces factor as the product of two rings in a way that forces
the principally polarized surfaces to split as a product of polarized elliptic
curves.
(See~\cite[Thm.~3.4]{maisner-nart}.)
\end{enumerate}

Our proof of Theorem~\ref{T-main} uses polarizations, but in
contrast to the last two methods
it is in essence a counting argument.  There are
many principally polarized abelian surfaces with characteristic polynomial $f_q$;
we compute exactly how many there are, and we compare this number to the number of
geometrically split principally polarized abelian surfaces with characteristic
polynomial $f_q$ that we obtain from a simple construction.
The first number turns out to be
equal to the second, so we conclude that every principally polarized
abelian surface with characteristic polynomial $f_q$ splits over the algebraic
closure and is hence not a Jacobian.

Lenstra, Pila, and Pomerance~\cite{lenstra-pila-pomerance}
give lower bounds for the number of principally polarized abelian
surfaces with characteristic polynomials of a certain form.  Our argument is
closely related to theirs, and in fact almost all of the ingredients of
our proof of Proposition~\ref{P-PPAS} (below) appear in their paper.  
The only reason we get an exact formula for the number of
principally polarized surfaces with characteristic polynomial~$f_q$, rather
than a lower bound along the lines of~\cite[Prop.~8.2]{lenstra-pila-pomerance},
is that every principally polarized variety we must consider
has an endomorphism ring that is Gorenstein.  We recommend the 
paper~\cite{lenstra-pila-pomerance} to those readers who would like
to see a more general application of the arguments we give
in~Section~\ref{S-PPAS}.

Maisner and Nart \cite[\S 3]{maisner-nart}
raise the question of determining which isogeny 
classes of simple non-supersingular abelian surfaces contain Jacobians.
They note that the only such isogeny classes for which the question
is presently unresolved have characteristic polynomials
of the form $f = x^4 + (a - 2q)x^2 + q^2$ with $0<a<4q$,
so it is natural to ask whether
our argument can be applied to any other polynomials of this form.
One can show that our argument also works for the case $a=1$, but in
that case there is a much simpler explanation for why the
corresponding $f$ is not the characteristic polynomial of a 
curve (see~\cite{howe:em}).  For every other allowable value of~$a$
our argument fails, because it depends (somewhat implicitly)
on the ring $\BZ[\sqrt{-a}]$ being the full ring of integers of
a field of class number one, and this happens only when $a=1$ or $a=2$.
On the other hand, it should be possible to adapt our argument to show that
except for these two cases,
the characteristic polynomials given above {\em do\/} come from curves.

In Section~\ref{S-outline} we give a detailed outline of our proof.
In Section~\ref{S-norm} we prove a technical result about the existence
of a well-behaved norm map between two class groups that will be important
in our argument.
In Section~\ref{S-PPAS} we count the number of principally polarized
abelian surfaces that have characteristic polynomial~$f_q$, and
in Section~\ref{S-split} we count the number of these polarized
surfaces that split over the quadratic extension of the base field.
Finally,
in Section~\ref{S-Brauer} we show that the two counts are equal, which 
shows that there is no curve with characteristic polynomial~$f_q$.

\begin{ack}
The author thanks Hendrik Lenstra for helpful conversations about 
the techniques used in this paper.
\end{ack}

\section{Outline of the proof}
\label{S-outline} 

Suppose $D$ is a positive squarefree integer that is congruent to $1$ modulo~$4$.
Let $K$ be the number field $\BQ(\sqrt{-2}, \sqrt{-D})$, let $\Kp$
be the maximal real subfield $\BQ(\sqrt{2D})$ of $K$, and let $L$ be the
imaginary quadratic subfield $\BQ(-D)$ of $K$.  In our notation for these fields
we have intentionally suppressed the dependence on~$D$, simply to keep the notation
from getting out of hand.
The number field $K$ and its subfields
will be important to our argument because if $A$ is an abelian surface over $\BF_q$
with characteristic polynomial $f_q$, and if one writes $2q-1 = F^2 D$ with $D$
squarefree, then $(\End A)\otimes\BQ \cong K$.  In fact, the isomorphism
can be chosen so that the Frobenius endomorphism on $A$ corresponds to
the element
$(F\sqrt{2D} + \sqrt{-2})/2$
of $K$.

Let $\CO$ denote the ring of integers $\BZ[\sqrt{-2}]$ of $\BQ(\sqrt{-2})$ and 
let $w$ denote the element $(\sqrt{2D} + \sqrt{-2})/2$ of $K$.  For every 
odd positive integer $f$ let $R_f$ be the subring $\CO + fw\CO$ of $K$,
let $\Rp_f$ be the subring $\BZ[f\sqrt{2D}]$ of $\Kp$, and let
$S_f$ be the subring $\BZ[f\sqrt{-D}]$ of $L$. 
Note that when $f=1$ these rings are the maximal orders of their quotient fields.
Let $U_f$ denote the group of units of $R_f$ and let $\Up_f$ denote the
group of totally positive units of $\Rp_f$.
(Recall that an element $x$ of a number field $M$ is {\em totally positive\/}
if $\varphi(x)>0$ for all embeddings $\varphi\colon M\to\BR$.)

The {\em class group} of an order $R$, denoted $\Cl R$, is the quotient
of the group of invertible fractional ideals of $R$ by the subgroup of 
principal fractional ideals.  The {\em narrow class group} of $R$, denoted
$\Clp(R)$, is the quotient of the group of invertible fractional ideals of $R$
by the subgroup of principal ideals that can be generated by a
totally positive element of the quotient field of $R$.

\begin{proposition}
\label{P-PPAS}
Let $q$ be a power of an odd prime and write $2q-1 = F^2 D$, where $D$ is squarefree.
The number of principally polarized abelian surfaces $(A,\lambda)$ over~$\Fq$
{\rm(}up to isomorphism over $\BF_q${\rm)}
such that $A$ has characteristic polynomial $f_q$ is equal to the sum
$$\sum_{f\mid F} [\Up_f : N(U_f)]\frac{\#\Cl R_f}{\#\Clp \Rp_f},$$
where $N$ denotes the norm map from $K$ to $\Kp$.
\end{proposition}

We see that there do exist principally polarized abelian surfaces with characteristic
polynomial $f_q$.  The only way that these polarized surfaces could fail to be 
Jacobians is if they split geometrically as products of polarized elliptic curves.
Thus, 
in order to prove Theorem~\ref{T-main} we must come up with a large enough supply of
these split surfaces to account for the surfaces enumerated in Proposition~\ref{P-PPAS}.
Maisner and Nart note that every abelian surface over $\BF_q$
with characteristic polynomial $f_q$ splits up to isogeny over $\BF_{q^2}$, so
it is natural to try to construct the requisite split surfaces using elliptic
curves over $\BF_{q^2}$.  (In fact, it is {\em necessary\/} to use elliptic curves
over $\BF_{q^2}$, as is shown by Proposition~\ref{P-splitting} below.)

Let $g_q$ be the polynomial $x^2 + (2 - 2q)x + q^2$, so that $f_q(x) = g_q(x^2)$.
The polynomial $g_q$ is the characteristic polynomial of an isogeny class $\CC$ of
elliptic curves over $\BF_{q^2}$.  Suppose $E$ is an elliptic curve
in this isogeny class and let $E^{(q)}$ be its Galois conjugate over $\BF_q$.
Then the product surface $A = E \times E^{(q)}$ has a natural principal polarization
$\lambda$, namely the product polarization.  It is easy to check that the obvious
isomorphism $f: (A,\lambda) \to (A,\lambda)^{(q)}$ satisfies $f^{(q)} \circ f = 1$,
so we can descend $(A,\lambda)$ to $\BF_q$.  Let us denote this geometrically split
principally polarized surface over $\BF_q$ by $S(E)$.
We find that the characteristic polynomial of Frobenius of $S(E)$ is~$f_q$,
and it is easy to check that $S(E_1)\cong S(E_2)$ as polarized varieties
if and only if either $E_2\cong E_1$ or $E_2\cong E_1^{(q)}$.  
Thus, if the action of $\Gal(\BF_{q^2}/\BF_q)$ on the 
isogeny class $\CC$ breaks $\CC$ into $n_1$ orbits of size $1$ and $n_2$ orbits 
of size $2$, we construct $n_1 + n_2$ distinct split surfaces
with characteristic polynomial~$f_q$ in this manner.

\begin{proposition}
\label{P-EC}
Let $q$ be a power of an odd prime and write $2q-1 = F^2 D$, where $D$ is squarefree.
The number of elliptic curves over $\BF_{q^2}$ with characteristic polynomial $g_q$
is equal to the sum
$$\sum_{f\mid F} \#\Cl S_f .$$
None of these curves is isomorphic to its Galois conjugate over $\BF_q$, unless $D=1$,
in which case exactly one of the curves is isomorphic to its conjugate.
\end{proposition}

Let $\eps$ denote the function on the integers such that $\eps(x) = 1$ if $x = 1$
and $\eps(x) = 0$ otherwise.
Proposition~\ref{P-EC} shows that the number of split surfaces obtained from 
$\CC$ is equal to
$$\frac{1}{2} \left( \eps(D) + \sum_{f\mid F} \#\Cl S_f \right)
=
\sum_{f\mid F}\frac{\eps(fD) + \#\Cl S_f}{2}.$$

We complete the proof of Theorem~\ref{T-main} by proving a purely number-theoretic
result about class groups, based on the Brauer class-number relations for the
extension~$K/\BQ$.

\begin{proposition}
\label{P-Brauer}
Let $D$ be a positive squarefree integer that is congruent to $1$ modulo $4$
and let $f$ be an odd positive integer.  Then 
$$[\Up_f:N(U_f)] \frac{\#\Cl R_f}{\#\Clp \Rp_f} = \frac{\eps(fD) + \#\Cl S_f}{2}.$$
\end{proposition}

We will prove Propositions~\ref{P-PPAS}, \ref{P-EC}, and~\ref{P-Brauer} in the
following sections.  Together with the discussion above, they provide a proof
of Theorem~\ref{T-main}.

\section{The norm map on class groups}
\label{S-norm}

In the course of our arguments we will need to know that there is a norm
map from the class group of $R_f$ to the narrow class group of $\Rp_f$,
and that this norm map is surjective.  We will prove these facts in this section.

If $T$ is an order in a number field,  we let $I(T)$ denote the group of 
invertible fractional $T$-ideals.
Note that $R_f = \Rp_f[(\sqrt{-2}+f\sqrt{2D})/2]$ so that 
$R_f$ is flat over~$\Rp_f$;
it follows that the tensoring-with-$R_f$ map from $I(\Rp_f)$ to $I(R_f)$ 
is an injection.
If $\Aid$ is an ideal of $R_f$ we let $\Aidbar$ denote its complex conjugate.

\begin{lemma}
\label{L-norm}
Suppose $\Aid\in I(R_f).$ Then $\Aid\Aidbar\in I(R_f)$ lies in the image of $I(\Rp_f)$.
\end{lemma}

\begin{proof}
Let $\Tor I(R_f)$ denote the torsion subgroup of $I(R_f)$.
A result of Dedekind~\cite{dedekind}, as
reinterpreted by Sands~\cite{sands}, says that the
tensoring-with-$R_1$ map $I(R_f)\to I(R_1)$ is surjective
and has kernel~$\Tor I(R_f)$.  In other words, 
there is an exact sequence
$$0 \to \Tor I(R_f) \to I(R_f) \to I(R_1) \to 0,$$
where the second map is inclusion and the third is tensoring-with-$R_1$.
Likewise, we have a sequence
$$0 \to \Tor I(\Rp_f) \to I(\Rp_f) \to I(\Rp_1) \to 0.$$
Let $\Bid$ be the image of $\Aid$ in $I(R_1)$.  Then the usual theory
of ideals in number fields shows that there is an ideal $\Cid\in I(\Rp_1)$
such that the image of $\Cid$ in $I(R_1)$ is equal to $\Bid\Bidbar$.
Let $\Did$ be an element of $I(\Rp_f)$ that maps to $\Cid$.
Let $\Eid = \Aid\Aidbar\Did^{-1}\in I(R_f)$.  Then $\Eid$
maps to the identity of $I(R_1)$, and hence lies in $\Tor I(R_f)$. 
We will be done if we can show that $\Eid$ lies in the image of 
$\Tor I(\Rp_f)$.

Let $\fid = f R_1$.  Note that $\fid$ is an ideal of $R_1$ and of $R_f$.
Sands shows that 
$$\Tor I(R_f) \cong (R_1\bmod \fid)^* / (R_f\bmod \fid)^*.$$
Similarly, if we let $\fidp = f \Rp_1$, then
$$\Tor I(\Rp_f) \cong (\Rp_1\bmod \fidp)^* / (\Rp_f\bmod \fidp)^*.$$
If we apply Galois cohomology to the sequence
$$0\to (R_f\bmod \fid)^*\to(R_1\bmod \fid)^*\to\Tor I(R_f)\to 0$$
and note that 
$$H^1(\Gal(K/\Kp), (R_f\bmod \fid)^*) = 
  H^1(\Gal(\BQ(\sqrt{-2})/\BQ, (\CO\bmod f\CO)^*) = 0,$$ 
we see that every element of $\Tor I(R_f)$ 
fixed by conjugation is represented by an element of $(R_1\bmod \fid)^*$
that is fixed by conjugation.  In particular, $\Eid$ is represented by
an element of $(R_1\bmod \fid)^*$ that is fixed by conjugation.  But
every such element comes from $(\Rp_1\bmod\fidp)^*$ because $f$ is odd.
Thus, $\Eid$ lies in the image of $\Tor I(\Rp_f)$, and we are done.
\end{proof}

Lemma~\ref{L-norm} shows that the map $\Aid\mapsto\Aid\Aidbar$ can be
viewed as a map from $I(R_f)$ to $I(\Rp_f)$.  Piecing this map together 
with the norm maps on $I(R_1)$ and $(R_1 \bmod \fid)^*$, we obtain
a diagram
$$\begin{matrix}
0 & \lra &  \Tor I(R_f)  & \lra &  I(R_f)  & \lra &  I(R_1)  & \lra & 0\\
  &      &      \bda     &      &   \bda   &      &   \bda   &      &  \\
0 & \lra & \Tor I(\Rp_f) & \lra & I(\Rp_f) & \lra & I(\Rp_1) & \lra & 0\\
\end{matrix}$$
This last diagram gives us the following diagram of class groups, in which
the vertical arrows are norm maps:
\begin{equation}
\label{EQ-class-group-diagram}
\begin{matrix}
0 & \lra &   D_f & \lra &   \Cl R_f  & \lra &   \Cl R_1  & \lra & 0\\
  &      &  \bda &      &    \bda    &      &    \bda    &      &  \\
0 & \lra & \Dp_f & \lra & \Clp \Rp_f & \lra & \Clp \Rp_1 & \lra & 0\\
\end{matrix}
\end{equation}
Here we have set
\begin{align*}
D_f & = \Tor I(R_f) / (U_1\bmod\fid) \\
    & = (R_1\bmod \fid)^* / ((R_f\bmod \fid)^*\cdot (U_1\bmod\fid))\\
\intertext{and}
\Dp_f & = \Tor I(\Rp_f) / (\Up_1\bmod\fidp)\\
      & = (\Rp_1\bmod \fidp)^* / ((\Rp_f\bmod \fidp)^* \cdot (\Up_1 \bmod \fidp)).
\end{align*}

We will have call to use the following basic fact about the norm from
$\Cl R_f$ to $\Clp \Rp_f$.
\begin{lemma}
\label{L-surjective-norm}
The norm map $\Cl R_f \to \Clp \Rp_f$ is surjective.
\end{lemma}

\begin{proof}
Note that the prime $2$ of $\BQ$ is totally ramified in the extension $K/\BQ$,
so the field extension $K/\Kp$ is ramified at a finite prime.
It follows from class field theory that $\Cl R_1 \to \Clp \Rp_1$ is surjective;
see~\cite[Prop.~10.1]{howe:tams}, for example.
Next note that the norm map from $(R_1\bmod \fid)^*$ to $(\Rp_1\bmod \fidp)^*$
is surjective, so that the norm from $D_f$ to $\Dp_f$ is surjective also.
Then we see from diagram~(\ref{EQ-class-group-diagram}) that 
the norm map $\Cl R_1 \to \Clp \Rp_1$ is surjective as well.
\end{proof}

\begin{rem}
The existence of a norm map $\Cl R\to\Clp\Rp$ is proven at the beginning 
of \S6 of~\cite{lenstra-pila-pomerance} through an identification of
the class groups with certain quotients of profinite groups.
\end{rem}

\section{The number of principally polarized surfaces}
\label{S-PPAS}

We will make heavy use of Deligne's equivalence of categories (see~\cite{deligne})
between the category of ordinary abelian varieties over $\BF_q$ and a certain
category of modules, called {\em Deligne modules\/} in~\cite{howe:tams}.
Deligne's category equivalence was expanded upon in~\cite{howe:tams}, so that
several geometric notions were translated into the category of Deligne modules.
In particular, there is a notion of a polarization of a Deligne module.  We will
count principally polarized abelian surfaces in the isogeny class determined by $f_q$
by counting principally polarized Deligne modules in the isogeny class
determined by $f_q$.  

We begin by reviewing the concepts of~\cite{howe:tams}.
However, we will simplify matters by restricting our discussion to Deligne modules
associated to irreducible polynomials.  Thus, the definitions we give below will
not be identical to the ones found in~\cite{howe:tams}, but will be equivalent
to the original definitions in our special case.

Let $q$ be a power of a prime $p$.
An {\em ordinary Weil $q$-polynomial\/} 
is a monic polynomial in $\BZ[x]$ of even degree
whose middle coefficient is coprime to $q$ and all of whose roots in the complex
numbers have magnitude $\sqrt{q}$.  (For example, $f_q$ is an ordinary Weil 
$q$-polynomial, and $g_q$ is an ordinary Weil $q^2$-polynomial.)
Suppose that $f$ is an irreducible
ordinary Weil $q$-polynomial.  Let $K$ be the number field $\BQ[x]/(f)$,
let $\pi$ be the image of $x$ in $K$, and let $\pibar = q/\pi$.
The number field $K$ is a totally imaginary quadratic extension of a totally real
field $\Kp$, and $\pibar$ is the complex conjugate of $\pi$.
Let $R$ be the order $\BZ[\pi,\pibar]$ of $K$.  A {\em Deligne module\/}
associated to $f$ is a finitely-generated non-zero sub-$R$-module of $K$.
A {\em morphism\/} from a Deligne module $\Aid$ to a Deligne module $\Bid$
is an element $\alpha$ of $K$ such that $\alpha\Aid\subseteq\Bid$.
An {\em isogeny\/} is a nonzero morphism.  The {\em degree\/} of an isogeny 
$\alpha:\Aid\to\Bid$ is the index of $\alpha\Aid$ in $\Bid$.

The {\em dual\/} $\hat{\Aid}$ of a Deligne module $\Aid$ is the complex conjugate of 
the dual of $\Aid$ under the trace pairing.  Suppose we pick, once and for all,
a $p$-adic  valuation $\nu$ on the algebraic closure of $\BQ$ in $\BC$.  Choose a
totally imaginary element $\iota$ of $K$ with the property that for every embedding
$\varphi:K\to\BC$ with $\nu(\varphi(\pi)) > 0$ the real number $\varphi(\iota) / i$
is positive.  A {\em polarization\/} of a
Deligne module $\Aid$ is an isogeny from $\Aid$ to $\hat{\Aid}$ of the
form $\iota\alpha$, where $\alpha$ is a totally positive element of the field $\Kp$.
A polarization is {\em principal\/} if it is an isomorphism.

Part of the main result of~\cite[\S 4]{howe:tams} is that Deligne's equivalence of
categories from the category of abelian varieties over $\Fq$ with characteristic
polynomial $f$ to the category of Deligne modules associated to $f$ takes isogenies
to isogenies, dual varieties to dual modules, and polarizations to polarizations.
Thus, to count principally polarized varieties in the isogeny class determined
by~$f$, it is sufficient to count principally polarized Deligne modules associated
to~$f$.

So suppose $\Aid$ is a principally polarized Deligne module associated to $f_q$.
Let $A = \End\Aid$.  Since a lattice and its trace dual have the same 
endomorphism ring, we see that $\End\hat{\Aid}$ is the complex conjugate of
$\End \Aid$.
Since $\Aid$ is supposed to be isomorphic to $\hat{\Aid}$, we see that $A$ is stable
under complex conjugation.  We also know that $A\supseteq R = \BZ[\pi,\pibar]$, 
and since $\pi$ can be taken to be $(F\sqrt{2D} + \sqrt{-2})/2$ we see that
$A$ contains $\pi - \pibar = \sqrt{-2}$.  Thus $A$ contains $\CO$.
It follows that $A = \CO \oplus w I$ for some ideal $I$ of $\CO$, and since
$A$ is stable under complex conjugation we must have $I = \overline{I}$.
Thus, either $I = f\CO$ for an integer $f$, or $I = f\sqrt{-2}\CO$ for an
integer $f$.  But $A$ contains $\pi = (1 - F)\sqrt{-2} / 2 + F w$, and since
$F$ is odd we see that we cannot have $I = f\sqrt{-2}\CO$.  It follows that $I = f\CO$
for some divisor of $F$; in other words, $A = R_f$ for some divisor $f$ of $F$.

First we count the number of Deligne modules with endomorphism ring
equal to~$R_f$.

\begin{lemma}
\label{L-SwithRf}
Let $f$ be a divisor of $F$.
The number of isomorphism classes of
Deligne modules $\Aid$ with $\End\Aid = R_f$ is equal to $\#\Cl R_f$.
\end{lemma}

\begin{proof}
First we note that $R_f$ is a Gorenstein ring because it is 
a complete intersection over $\BZ$.  (In fact, $R_f$ is the quotient
of $\BZ[x,y]$ by the ideal generated by the regular sequence
$(x^2 + 2, y^2 - fxy - f^2(D+1)/2)$, the image of $x$ being $\sqrt{-2}$
and the image of $y$ being $fw$.)
Since $R_f$ is Gorenstein, every fractional $R_f$-ideal $\Aid$
with $\End\Aid = R_f$ is actually an invertible $R_f$-ideal.
Thus, the number of isomorphism classes of
Deligne modules $\Aid$ with $\End\Aid = R_f$ is equal to $\#\Cl R_f$.
\end{proof}

Next, we give a method for deciding which elements of $\Cl R_f$
give Deligne modules that admit principal polarizations. 

\begin{lemma}
\label{L-PPASwithRf}
Let $f$ be a divisor of $F$.
A Deligne module $\Aid$ with $\End\Aid = R_f$ has a principal
polarization if and only if its class in $\Cl R_f$ lies in the kernel of the
norm map from $\Cl R_f$ to $\Cl^+ \Rp_f$.
The number of isomorphism classes of 
Deligne modules $\Aid$ with $\End\Aid = R_f$ that have a principal
polarization is equal to the quotient
$$\frac{\#\Cl R_f}{\#\Cl^+ \Rp_f}.$$
\end{lemma}

\begin{proof}
Before we discuss polarizations, we must pick an appropriate totally imaginary
element $\iota$ of $K$, which requires that we pick a $p$-adic valuation $\nu$
on the algebraic closure of $\BQ$ in~$\BC$.  To pick $\nu$, we pick our
favorite embedding $\varphi$ of $K$ into~$\BC$, pick a prime $\pid$ of $K$
containing $\pi$, and choose $\nu$ to be any valuation that extends the 
$\pid$-adic valuation on $\varphi(K)$.  
Now the condition on $\iota$ is that $\varphi(\iota)/i$ should be positive
and that $\varphi(\iota^\sigma / \iota)$ should be positive for every
$\sigma\in\Gal(K/\BQ)$ such that~$\pi^\sigma\in\pid$.
Since $\pi^2\in \BQ(\sqrt{-D})$, these 
$\sigma$ are precisely the automorphisms that fix $\BQ(\sqrt{-D})$.
Thus we can pick $\iota = \pm1/(4\sqrt{-D})$, with the sign chosen
so that $\varphi(\iota)/i$ is positive.

Suppose $\Aid$ is a Deligne module with $\End \Aid = R_f$.  How can we tell
whether $\Aid$ has a principal polarization?  First we note that the
different of the order $R_f$ is $4 f \sqrt{-D}$, so the trace dual
of $\Aid$ is $\Aid^{-1} / (4 f \sqrt{-D})$.  It follows that
$\hat{\Aid} = \Aidbar^{-1} / (4 f \sqrt{-D})$.   Then
$\Aid$ has a principal polarization if and only if there is a totally
positive element $\alpha$ of $\Kp$ such that
$$\alpha\iota \Aid = \hat{\Aid} = \Aidbar^{-1} / (4 f \sqrt{-D}),$$
which is equivalent to 
$$\Aid\Aidbar = \left(\frac{1}{f\alpha}\right) R_f.$$
But as the results of Section~\ref{S-norm} show, this last condition is 
equivalent to the condition that the class of $\Aid$ in $\Cl R_f$ be in the
kernel of the norm map from $\Cl R_f$ to $\Clp \Rp_f$.
This proves the first statement of the lemma.

The second statement of the lemma follows immediately from the first
because Lemma~\ref{L-surjective-norm} says that
the norm from $\Cl R_f$ to $\Cl^+ \Rp_f$ is surjective.
\end{proof}

Finally, we count how many principal polarizations a Deligne module has,
given that it has at least one.

\begin{lemma}
\label{L-CountPols}
Let $f$ be a divisor of $F$, and suppose $\Aid$ is a Deligne module
with $\End\Aid = R_f$ that has a principal polarization.  Then the 
number of non-isomorphic principal polarizations on $\Aid$
is equal to $[\Up_f : N(U_f)]$, where $N\colon U_f\to\Up_f$
is the norm map from the units of $R_f$ to the totally positive
units of $\Rp_f$.
\end{lemma}

\begin{proof}
Let $\lambda$ be a principal polarization of~$\Aid$.  Then the map
$u\mapsto u\lambda$ gives a bijection between the totally positive units
of $\Rp_f$ and the principal polarizations of~$\Aid$.  On the other hand,
if $\mu_1$ and $\mu_2$ are principal polarizations of $\Aid$,
then $(\Aid,\mu_1)$ and $(\Aid,\mu_2)$ are isomorphic polarized Deligne modules
if and only if there is a unit $u$ of $R_f$ such that $\mu_1 = u\overline{u}\mu_2$.
Thus, the isomorphism classes of principal polarizations of $\Aid$ correspond
to the elements of the group $\Up_f / N(U_f).$
\end{proof}

Together, Lemmas~\ref{L-SwithRf}, \ref{L-PPASwithRf}, and~\ref{L-CountPols}
prove Proposition~\ref{P-PPAS}.\qed

\section{The number of surfaces coming from elliptic curves}
\label{S-split}

In this section we count the number of elliptic curves over $\BF_{q^2}$ in the
isogeny class $\CC$ determined by $g_q = x^2 + (2 - 2q)x + q^2$, and we determine the
number of fixed points of the action of $\Gal \BF_{q^2} / \BF_q$ on $\CC$.

Let $\rho$ be a root of $g_q$ in the field $L = \BQ(\sqrt{-D})$ discussed above.
Note that $\BZ[\rho] = S_F$.  By a classical result of Deuring~\cite{deuring},
the number of elliptic curves in $\CC$ is given by the sum
$$\sum_{f\mid F} \#\Cl S_f,$$
and this is the main part of the statement of Proposition~\ref{P-EC}.  We
are left to prove the statement about the Galois action on $\CC$.

Suppose $E$ is an elliptic curve in $\CC$ that is isomorphic to its Galois
conjugate $E^{(q)}$.  Let $f\colon E \to E^{(q)}$ be an isomorphism, and
consider the automorphism $g = f^{(q)}\circ f$ of $E$.

Since $E$ is ordinary, the automorphism group of $E$ is cyclic of order $2$, $4$, 
or $6$.  Order $6$ is impossible, because then $L$ would have to be $\BQ(\sqrt{-3})$,
whereas we know that $D$ is congruent to $1$ modulo $4$.   If $E$ has an automorphism
of order $4$ then $D = 1$ and $E$ has $j$-invariant $0$; let us consider this case
later, and focus for now on the case where $\Aut E =\{\pm 1\}$.

If $g = 1$ then we can descend $E$ to $\BF_q$; that is, we find that there is an
elliptic curve $F$ over $\BF_q$ that becomes isomorphic to $E$ when we extend
the base field to~$\BF_{q^2}$.  If $g = -1$ then we find something slightly
different: there is an elliptic curve $F$ over $\BF_q$ that becomes isomorphic
to $E$ when we extend the base field to~$\BF_{q^4}$, or in other words, when we
extend the base field to $\BF_{q^2}$ the curve $F$ becomes isomorphic to the
quadratic twist of $E$.

Suppose the characteristic polynomial of $F$ over $\BF_q$ is $x^2 - sx + q$.  Then 
the characteristic polynomial of $F$ over $\BF_{q^2}$ is 
$x^2 + (2q - s^2) x + q^2$.  Since $F$ becomes isomorphic over $\BF_{q^2}$
to either $E$ or its quadratic twist, this last polynomial must be either
$g_q(x)$ or $g_q(-x)$.  Thus we must have either $2q  - s^2 = 2 - 2q$ or
$2q - s^2 = 2q - 2$.  The latter condition requires $s^2 = 2$, which is
clearly impossible.  The former condition requires that $4q = s^2 + 2$,
which cannot hold because it is impossible modulo $4$.  Thus, we obtain
a contradiction if $\Aut E = \{\pm1\}.$

If $D = 1$ then the isogeny class $\CC$ does contain an elliptic
curve $E$ with $\End E = \BZ[i]$, and this curve is unique because the
class number of $\BZ[i]$ is $1$.  Thus, this $E$ is isomorphic to its 
Galois conjugate.

This proves Proposition~\ref{P-EC}. \qed

While we do not need the following result for our argument, it is good
to point out that the geometrically split polarized surfaces we counted 
above are the only possible split polarized surfaces in the isogeny class
defined by $f_q$.  

\begin{proposition}
\label{P-splitting}
Suppose $(A,\lambda)$ is a principally polarized abelian surface over
a finite field $k$.  If $(A,\lambda)$ is not the polarized Jacobian of
a genus-$2$ curve over $k$, then over the quadratic extension of $k$
the polarized surface $(A,\lambda)$ may be written as a product of
two principally polarized elliptic curves.
\end{proposition}

\begin{proof}
We know that over {\em some\/} finite extension $\ell$ of $k$, say of degree $n$,
the polarized surface $(A,\lambda)$ can be
written as a product of polarized elliptic curves $(E_1,\mu_1)$
and $(E_2,\mu_2)$, together with descent data.
If we let $\sigma$ be the Frobenius automorphism of the extension $\ell/k$,
the descent data is an isomorphism 
$f :E_1\times E_2 \to E_1^\sigma\times E_2^\sigma$
that respects the polarizations $\mu_1\times\mu_2$ and $\mu_1^\sigma\times\mu_2^\sigma$
and that satisfies the relation 
$$ 1 = f^{\sigma^{n-1}} \circ \cdots\circ f^\sigma\circ f.$$
Using the fact that $f$ respects the polarizations one can show that $f$
can be represented by a $2\times 2$ matrix of isogenies that is either 
diagonal or anti-diagonal.  In the former case, the descent data gives rise to
descent data for $(E_1,\mu_1)$ and $(E_2,\mu_2)$ that shows that $(A,\lambda)$
splits already over $k$.  In the latter case, we see that $n$ must be even, 
and that $f^\sigma\circ f$ gives descent data for $(A,\lambda)$ from $\ell$ to
the quadratic extension of~$k$.  Again we find that the descent data for the
surface gives descent
data for the polarized elliptic curves, so $(A,\lambda)$ splits over the
quadratic extension of~$k$.
\end{proof}

Proposition~\ref{P-splitting} is almost the same as
Theorem~3.4 from~\cite{maisner-nart},
but that theorem speaks only of the splitting of $A$, not of $(A,\lambda)$.

\section{The Brauer relations}
\label{S-Brauer}

In this section we will prove Proposition~\ref{P-Brauer}. 
Our proof is based on two commutative diagrams.  The first diagram is
simply diagram~(\ref{EQ-class-group-diagram}) with the kernels $C$, $B$,
and $A$ of the vertical maps added in.
\begin{equation}
\label{EQ-big-class-group-diagram}
\begin{matrix}
  &      &   0   &      &      0     &      &      0     &      &  \\
  &      & \bda  &      &    \bda    &      &    \bda    &      &  \\
0 & \lra &   C   & \lra &      B     & \lra &      A     & \lra & 0\\
  &      & \bda  &      &    \bda    &      &    \bda    &      &  \\
0 & \lra &  D_f  & \lra &  \Cl R_f   & \lra &  \Cl R_1   & \lra & 0\\
  &      & \bda  &      &    \bda    &      &    \bda    &      &  \\
0 & \lra & \Dp_f & \lra & \Clp \Rp_f & \lra & \Clp \Rp_1 & \lra & 0\\
  &      & \bda  &      &    \bda    &      &    \bda    &      &  \\
  &      &   0   &      &      0     &      &      0     &      &  \\
\end{matrix}
\end{equation}
The second diagram comes from the definitions of the groups in the
left-most column of diagram~(\ref{EQ-big-class-group-diagram}).
\begin{equation}
\label{EQ-big-order-diagram}
\begin{matrix}
  &      &   0   &      &      0     &      &      0     &      &  \\
  &      & \bda  &      &    \bda    &      &    \bda    &      &  \\
0 & \lra &   F   & \lra &      E     & \lra &      C     &      &  \\
  &      & \bda  &      &    \bda    &      &    \bda    &      &  \\
0 & \lra &  G_f  & \lra & (R_1\bmod\fid)^*/(R_f\bmod\fid)^*
                                     & \lra &     D_f    & \lra & 0\\
  &      & \bda  &      &    \bda    &      &    \bda    &      &  \\
0 & \lra & \Gp_f & \lra & (\Rp_1\bmod\fidp)^*/(\Rp_f\bmod\fidp)^*
                                     & \lra &    \Dp_f   & \lra & 0\\
  &      & \bda  &      &    \bda    &      &    \bda    &      &  \\
  &      &   H   & \lra &      0     &      &      0     &      &  \\
  &      & \bda  &      &            &      &            &      &  \\
  &      &   0   &      &            &      &            &      &  \\
\end{matrix}
\end{equation}
Here we have set
\begin{align*}
G_f & = \frac{U_1\bmod\fid}{(U_1\bmod\fid)\cap(R_f\bmod\fid)^*}\\
\intertext{and}
\Gp_f &= \frac{\Up_1\bmod\fidp}{(\Up_1\bmod\fidp)\cap(\Rp_f\bmod\fidp)^*},
\end{align*}
we denote by $F$ and $H$ the kernel and cokernel of the map $G_f\ra\Gp_f$,
and we denote by $E$ the kernel of the middle vertical map.

To prove Proposition~\ref{P-Brauer} we would like to calculate $[\Up_f:N(U_f)]\#B$.  
Diagram~(\ref{EQ-big-class-group-diagram}) shows that $\#B = \#A\#C$,
and diagram~(\ref{EQ-big-order-diagram}) shows that 
$\#C = \#E\#H/\#F$ and that $\#H/\#F = \#\Gp_f/\#G_f$.
Thus, $\#B = \#A\#E\#\Gp_f/\#G_f$,
and we have
\begin{equation}
\label{EQ-main-equation}
[\Up_f:N(U_f)]\#B = \big([\Up_1:N(U_1)]\#A\big)
                    \big(\#E\big)
                    \frac{[N(U_1) : N(U_f)]\#\Gp_f}
                         {[\Up_1 : \Up_f]\#G_f}.
\end{equation}
We will evaluate the right-hand side in three steps.  The steps are
slightly different when $D=1$, so for now we will assume that $D>1$.

\begin{lemma}
\label{L-level-one}
If $D>1$ then $[\Up_1:N(U_1)]\#A = (\#\Cl S_1) /2.$
\end{lemma}

\begin{proof}
It is easy to check that when $D>1$ the only roots of unity in $K$ are $\pm1$.
Let $u$ be a fundamental unit of $\Kp$ that is positive under at least one
of the embeddings of $\Kp$ into $\BR$, and let $v$ be a fundamental unit of $K$.
We take $v = u$ if the unit groups of $K$ and $\Kp$ are identical. 
Note that if the unit group of $K$ strictly contains the unit group of $\Kp$
then $N(v) = u$ and $u$ is totally positive.

We will say that we are in Case~1 if $v = u$ and $u$ is totally positive;
in Case~2 if $v=u$ and $u$ is not totally positive; and in Case~3
if $v\neq u$.  Note that in Cases 1 and 2 the regulator of $K$ is twice that
of $\Kp$, while in Case 3 the two regulators are equal.

For every number field $M$ let $q(M)$ denote the product of the class number
and the regulator of $M$, divided by the number of roots of unity in $M$.
The Brauer class number relations (\cite[\S VIII.7]{frohlich-taylor}),
applied to the biquadratic extension $K$ of $\BQ$, show that
$$ q(K) q(\BQ)^2 = q(\BQ(\sqrt{-2})) q(\Kp) q(L).$$
Since each of the fields contains exactly 2 roots of unity, and since
we know how the regulators of $K$ and $\Kp$ are related to one another,
we find that
$$\#\Cl R_1 = 
\begin{cases}
(\#\Cl \Rp_1)(\#\Cl S_1)/2 & \text{in Cases 1 and 2;}\\
(\#\Cl \Rp_1)(\#\Cl S_1) & \text{in Case 3.}\\
\end{cases}
$$
Now, in Cases 1 and 3 the narrow class number of $\Kp$ is twice the
class number of $\Kp$, while in Case~2 the narrow class number and the
class number are equal.  Thus we find that
$$\#\Cl R_1 = 
\begin{cases}
(\#\Clp \Rp_1)(\#\Cl S_1)/4 & \text{in Case 1;}\\
(\#\Clp \Rp_1)(\#\Cl S_1)/2 & \text{in Cases 2 and 3.}\\
\end{cases}
$$
Now we check that $[\Up_1 : N(U_1)]$ is $1$ in Cases 2 and 3, and is $2$ in Case~1.
Thus, 
$$[\Up_1 : N(U_1)] \#\Cl R_1 = (\#\Clp \Rp_1)(\#\Cl S_1)/2$$
in every case.
The lemma follows.
\end{proof}

\begin{lemma}
\label{L-unit-adjustment}
If $D>1$ then $\#G_f = [N(U_1) : N(U_f)]$ and $\#\Gp_f = [\Up_1 : \Up_f]$.
\end{lemma}

\begin{proof}
Recall that $U_1$ is generated by $-1$ and $v$. Since $-1\in R_f$, we see that
$\#G_f$ is equal to the smallest positive integer $i$ such that $v^i \bmod \fid$
lies in $R_f \bmod \fid$.  But since $\fid\subset R_f$, this $i$ is also the
smallest integer such $v^i\in R_f$, and this is clearly the index of
$U_f$ in $U_1$.  Since the norm kills only the elements $1$ and $-1$ in
$U_f$ and $U_1$, we see that $i$ is also the index of $N(U_f)$ in $N(U_1)$.

The proof of the second equality of the lemma is similar.
\end{proof}

\begin{lemma}
\label{L-class-number-adjustment}
If $D>1$ then $\#E = \#\Cl S_f / \#\Cl S_1.$
\end{lemma}

\begin{proof}
Since $S_1^*$ consists only of $\pm 1$, arguments as in 
Section~\ref{S-norm} show that
$$\frac{\#\Cl S_f}{\#\Cl S_1} = \frac{\#(S_1\bmod fS_1)^*}{\#(S_f\bmod fS_1)^*}.$$
On the other hand, the definition of $E$ shows that
$$\#E = \frac{ \#(R_1\bmod fR_1)^* \#(\Rp_f\bmod f\Rp_1)^*}
             { \#(R_f\bmod fR_1)^* \#(\Rp_1\bmod f\Rp_1)^*}.$$
If we note that
\begin{align*}
S_f\bmod fS_1 & \cong \BZ\bmod f\BZ\\
\Rp_f\bmod f\Rp_1 & \cong \BZ\bmod f\BZ\\
\intertext{and}
R_f \bmod fR_1 & \cong \CO\bmod f\CO,
\end{align*}
we see that what we want to show is that
$$\#(R_1/fR_1)^* (\#(\BZ/f\BZ)^*)^2 = \#(\CO/f\CO)^*\#(\Rp_1/f\Rp_1)^*\#(S_1/fS_1)^*$$
for every odd $f$.  It clearly suffices to check this when $f$ is a prime power.
We leave the details of this computation to the reader.
(But for example, we find that if $f=p^e$ for a prime
$p$ that splits completely in $K$, then the two sides of the equality above
are both equal to $(p-1)^6 p^{6e-6}$.)
\end{proof}

Combining Lemmas~\ref{L-level-one}, \ref{L-unit-adjustment}, and
\ref{L-class-number-adjustment} with equation~(\ref{EQ-main-equation}),
we find that 
$$[\Up_f:N(U_f)] \frac{\#\Cl R_f}{\#\Clp \Rp_f} = \#\Cl S_f/2$$
when $D>1$, as we were to show.

The argument when $D=1$ is entirely similar to the argument we have just given,
except that there are some complications because of the eighth roots of unity
in $K$ and the fourth roots of unity in $L$.  We leave the details of the
argument to the interested reader, but we will at least state the
required variants of Lemmas~\ref{L-level-one}, \ref{L-unit-adjustment}, and
\ref{L-class-number-adjustment}.

\begin{lemma}
\label{L-level-one-prime}
If $D=1$ then $[\Up_1:N(U_1)]\#A = \#\Cl S_1.$
\end{lemma}

\begin{lemma}
If $D=1$ and $f>1$ then $\#G_f = 4[N(U_1) : N(U_f)]$ and $\#\Gp_f = [\Up_1 : \Up_f]$.
\end{lemma}

\begin{lemma}
If $D=1$ and $f>1$ then $\#E = 2 \#\Cl S_f / \#\Cl S_1.$
\end{lemma}

When $D=1$ and $f=1$, Proposition~\ref{P-Brauer} is simply 
Lemma~\ref{L-level-one-prime}.  When $D=1$ and $f>1$, we can combine the
three lemmas above with equation~(\ref{EQ-main-equation}) to find that
$$[\Up_f:N(U_f)] \frac{\#\Cl R_f}{\#\Clp \Rp_f} = \#\Cl S_f/2.$$
This proves Proposition~\ref{P-Brauer}.
\qed

\end{document}